\documentclass[11pt, reqno]{amsart}
\usepackage[english]{babel}
\usepackage[table]{xcolor}
\usepackage[utf8]{inputenc}
\usepackage{amsmath, amsfonts, amsthm, amscd, amssymb, fullpage, fancyhdr, upgreek, amssymb, hyperref, enumerate, comment, siunitx, enumitem, graphicx, mathrsfs, mathtools, filecontents, xcolor,
microtype}
\usepackage{listings}
\usepackage[scr]{rsfso}
\usepackage{thmtools}
\usepackage{thm-restate}

\newcolumntype{P}[1]{>{\centering\arraybackslash}p{#1}}

\usepackage[parfill]{parskip}

\newtheorem{theorem}{Theorem}[section]

\newtheorem{lemma}[theorem]{Lemma}
\newtheorem{conjecture}[theorem]{Conjecture}

\theoremstyle{definition}
\newtheorem{definition}[theorem]{Definition}
\newtheorem{example}[theorem]{Example}
\newtheorem{remark}[theorem]{Remark}

\newcommand{\bburl}[1]{\textcolor{blue}{\url{#1}}}

\usepackage{tikz}
\usepackage{tkz-tab}
\usepackage{tkz-graph}
\usetikzlibrary{shapes.geometric,positioning}

\numberwithin{equation}{section}

\newcommand{\hr}[1]{\href{#1}{\url{#1}}}

\title{The Reversed Zeckendorf Game}

\author{Zo\"e X. Batterman}
\address{Pomona College}
\email{\textcolor{blue}{\href{mailto:zxba2020@mymail.pomona.edu}{zxba2020@mymail.pomona.edu}}}

\author{Aditya Jambhale}
\address{University of Cambridge}
\email{\textcolor{blue}{\href{mailto:aj644@cam.ac.uk}{aj644@cam.ac.uk}}}

\author{Steven J. Miller}
\address{Williams College}
\email{\textcolor{blue}{\href{mailto:sjm1@williams.edu}{sjm1@williams.edu}},  \textcolor{blue}{\href{Steven.Miller.MC.96@aya.yale.edu}{Steven.Miller.MC.96@aya.yale.edu}}}

\author{Akash L. Narayanan}
\address{Department of Mathematics, University of Michigan, Ann Arbor, MI 48104}
\email{\textcolor{blue}{\href{mailto:anaray@umich.edu}{anaray@umich.edu}}}

\author{Kishan Sharma}
\address{University of Cambridge}
\email{\textcolor{blue}{\href{mailto:kds43@cam.ac.uk}{kds43@cam.ac.uk}}}

\author{Andrew K. Yang}
\address{University of Cambridge}
\email{\textcolor{blue}{\href{mailto:aky30@cam.ac.uk}{aky30@cam.ac.uk}}}

\author{Chris Yao}
\address{Yale University}
\email{\textcolor{blue}{\href{mailto:chris.yao@yale.edu}{chris.yao@yale.edu}}}

\thanks{This work was completed during the 2023 SMALL REU program at Williams College. It was supported in part by NSF Grants DMS1561945 and DMS1659037, Williams College, and Churchill College, Cambridge.}


\date{\today}

\begin{document}

\maketitle

\begin{abstract} 
Zeckendorf \cite{Ze} proved that every natural number $n$ can be expressed uniquely as a sum of non-consecutive Fibonacci numbers, called its \emph{Zeckendorf decomposition}. Baird-Smith, Epstein, Flint, and Miller
\cite{original_paper} created the \emph{Zeckendorf game}, a two-player game played on partitions of $n$ into Fibonacci numbers which always terminates at a Zeckendorf decomposition, and proved that Player 2 has a winning strategy for $n\geq3$. Since their proof was non-constructive, other authors have studied the game to find a constructive winning strategy, and lacking success there turned to related problems. For example, Cheigh, Moura, Jeong, Duke, Milgrim, Miller, and Ngamlamai \cite{TowardsGauss} studied minimum and maximum game lengths and randomly played games.
We explore a new direction and introduce the \emph{reversed Zeckendorf game}, which starts at the ending state of the Zeckendorf game and flips all the moves, so the reversed game ends with all pieces in the first bin. We show that Player 1 has a winning strategy for $n = F_{i+1} + F_{i-2}$ and solve various modified games. 
\end{abstract}

\setcounter{tocdepth}{1}
\tableofcontents

\newpage
\section{Introduction and Main Results}
\subsection{History}
The Fibonacci numbers, which for uniqueness results in decompositions requires us to define them by $F_1 = 1$, $F_2 = 2$, and $F_{n+1} = F_n + F_{n-1}$, is a sequence with many interesting properties which have been widely studied. With this choice of indexing, Zeckendorf \cite{Ze} proved that every natural number $n$ has a unique \emph{Zeckendorf decomposition}, which expresses $n$ as the sum of distinct, non-adjacent Fibonacci numbers. Note the decomposition would no longer be unique if we considered there to be two ones or a zero in the Fibonacci sequence. An example of such a Zeckendorf decomposition is 
\begin{equation}
2024 \ = \ 1597 + 377 + 34 + 13 + 3.
\end{equation}
Building upon this, Baird-Smith, Epstein, Flint and Miller in \cite{original_paper} and \cite{BEFM2} created the \emph{Zeckendorf game}. The Zeckendorf game, which we will refer to as the \emph{forwards game}, starts with a natural number $n$. A game state consists of a partition of $n$ into Fibonacci numbers, constituents of which we call \emph{chips}. We say the collection of $F_i$'s in this partition is the $i$\textsuperscript{th} bin, with its cardinality $h_i$ called the \emph{height} of the bin. There are two types of moves.
\begin{enumerate}[label = (\roman*)]
\item \textit{Combine}: If $h_i>0$ and $h_{i-1}>0$, then the move is
\begin{align}
F_{i-1} + F_{i} &\ \mapsto\ F_{i+1}, \\
2F_1 &\ \mapsto\ F_2. \nonumber
\end{align}
In other words, we remove one chip from each of the $i$\textsuperscript{th} and $(i-1)$\textsuperscript{th} bins and add a chip to the $(i+1)$\textsuperscript{th} bin.
\item \textit{Split}: If $h_i\geq 2$ with $i> 2$, then we have the move
\begin{align}
2F_i &\ \mapsto\ F_{i-2} + F_{i+1}, \\
2F_2 &\ \mapsto\ F_{3} + F_{1}. \nonumber
\end{align}
\end{enumerate}
Note that each of these moves keeps the total sum of the values of the chips constant at $n$. The forwards game starts with $n$ ones and proceeds until a player can't make a move, in which case the player who can't move loses. This makes the Zeckendorf game a normal, impartial, combinatorial game.\footnote{The terms normal and impartial are defined in \textsection \ref{sec:reversed_games}.} 

The authors in \cite{original_paper} showed the game for $n$ always terminates at $n$'s Zeckendorf decomposition, and they showed, nonconstructively, Player 2 wins for all $n\geq 3$. After this result was published, many authors have investigated the properties of this game and various modifications of it such as possible game lengths, random games, etc; see \cite{CDHK+61}, \cite{CDHK+62}, \cite{BCD}, \cite{BEFM2}, \cite{TowardsGauss}, \cite{GMRVY}, \cite{LLMMSXZ}, and \cite{MSY}. 

We introduce another version of this game, the \emph{reversed Zeckendorf game}, where the players play with all the moves reversed and the starting and ending positions exchanged; for the full definition, see \textsection \ref{sec:reversed_zeck_game}. This game demonstrates a more complex winning structure than the forwards game.\footnote{Here, ``winning structure'' refers to when each player has a forced win.}

\subsection{Main Results}
First, in \textsection \ref{sec:reversed_zeck_game}, we analyze the reversed Zeckendorf game, look at its winning structure, and make further conjectures. As in the literature, we denote the golden ratio $\frac{1+\sqrt{5}}{2}$ by $\phi$. Moveover, for brevity in proofs, we will say that a player ``has a win'' (respectively ``has a loss'') by which we mean ``has a forced winning strategy'' (respectively ``the opponent has a force winning strategy''). We state our main results.
\begin{restatable}{theorem}{PlayerOneWinsInfinite}\label{thm:PlayerOneWinsInfinite}
Player 1 has a winning strategy for the reversed Zeckendorf game when 
\begin{align}
n = F_{i+1} + F_{i-2}
.\end{align}
\end{restatable}
We have both a nonconstructive and constructive proof for this theorem. The nonconstructive proof involves a strategy-stealing argument, and the constructive proof provides an explicit strategy.
\begin{restatable}{conjecture}{PlayerTwoWinsToo}\label{conj:infinitewin}
For the reversed Zeckendorf game, Player 2 has a winning strategy for infinitely many $n$.
\end{restatable}

\begin{restatable}{conjecture}{WinningPercentage}\label{conj:strong}
In the limit, the percent of the time Player 1 has a winning strategy for the reversed Zeckendorf game is $\varphi^{-1} \approx .618$.
\end{restatable}

Continuing in \textsection \ref{sec:otherfacts}, we transfer results from \cite{TowardsGauss} to get results about randomly played reversed Zeckendorf games as well as bounds on how long the game can last. These results are summarized here.

\begin{restatable}{theorem}{gamelengths}\label{thm:gamelengths}
Let $Z(n)$ be the number of terms in the Zeckendorf decomposition of $n$.
\begin{enumerate}[label = (\roman*)]
\item The shortest possible reversed Zeckendorf game is $n - Z(n)$.
\item An upper bound on the longest possible reversed Zeckendorf game is 
\begin{equation}
\lfloor \phi^2 n-Z_I(n) - 2Z(n)+\phi-1 \rfloor, 
\end{equation}
\end{enumerate}
where $Z_I(n)$ is the sum across the indices of the Fibonacci numbers in the Zeckendorf decomposition.
\end{restatable}
\begin{restatable}{theorem}{randomgames}\label{thm:randomgames}
    For any integer $Z \geq 1$ and $z \in \{0,1,\dots,Z-1\}$, we have
\begin{equation}
\lim_{N \to \infty} \mu_N(\text{game length equals } z \text{ mod } Z)\ =\ \lim_{N \to \infty} \mathbb{P }_N(\text{game length equals } z \text{ mod } Z)\ =\ \frac{1}{Z},
\end{equation}
where $\mu_N$ and $\mathbb{P}_N$ are two different probability measures, defined in \textsection \ref{sec:otherfacts}, being placed on the space of all possible games starting at the Zeckendorf decomposition of $N$. 
\end{restatable}
In \textsection \ref{sec:startpos}, we analyze winning strategies for different starting positions of the reversed Zeckendorf game. That is, we assume the starting position is not at the Zeckendorf decomposition but at some other partition of $n$ into Fibonacci numbers. Specifically, we completely solve the game when the starting position only involves ones, twos, and threes.
\begin{restatable}{theorem}{OneTwoThree}\label{thm:123}
    For any game starting with $a$ ones, $b$ twos, and $c$ threes, we have the following forced wins.
    \[
    \begin{tabular}{|c|c|c|c|c|}
    \hline
    a & b & c & & Player having forced win \\
    \hline
    Even & Even & Even & & Player 2 \\
    \hline
    Odd & Odd & Odd & & Player 1 \\
    \hline
    Even & Odd & Even & & Player 1 \\
    \hline 
    Odd & Even & Odd & & Player 1  \\
    \hline
    Odd & Even & Even & $a>c$& Player 2  \\
    \hline
    Odd & Even & Even &$a<c$ & Player 1  \\
    \hline
    Even & Even & Odd &  $a>c$& Player 1  \\
    \hline
    Even & Even & Odd &$a<c$ & Player 2  \\
    \hline
    Even & Odd & Odd & & Player 1  \\
    \hline
    Odd & Odd & Even & & Player 1  \\
    \hline
    \end{tabular}
    \]
\end{restatable}
Finally, in \textsection \ref{sec:buildup}, we conclude with a further modification of this game, which we call the \emph{build up game}, and we completely solve it. The proof involves Nim-like strategies---strategies that involve forcing the opponent along game states associated to a fixed residue class modulo some integer.

\section{Reversed Games}\label{sec:reversed_games}
Intuitively, two players playing a reversed game should look, if time is running backwards, like the forwards game. As such, we want a move in the reversed game from a state $A$ to a state $B$ to be legal if and only if moving from $B$ to $A$ is legal in the forwards game. To formalize this, we add the assumption that the game is \emph{normal and impartial}.  
\begin{definition}
We say a game is \emph{normal and impartial} when
\begin{enumerate}[label = (\roman*)]
\item the allowable moves depend only on the position and not on player order, and
\item the last player to move wins.
\end{enumerate}
\end{definition}
We can consider a normal, impartial game as a game played on a directed graph. Players take turns moving a game piece across the graph until a terminal node is reached (e.g., a node where there are no moves). 
\begin{definition}[\textbf{Reversed Game}]
Let $G$ be the associated directed graph of a normal, impartial game. Suppose that the graph has a unique loss node and unique starting node. 
Then the \emph{reverse} of that game is obtained by reversing all the edges in $G$ and playing from the loss node to the starting node. A player loses if they run out of moves.
\end{definition}\label{def:reverse}
\begin{remark}
Under this definition, the reverse of a game is also normal and impartial.
\end{remark}
\begin{example}
Consider the reversed game for standard Nim, where players take turns removing $1$, $2$, or $3$ from an natural number $n \equiv 1 \pmod 4$, never going below $1$, with the last player to move winning. Note standard Nim is a normal, impartial game with starting node $n$ and ending node $1$.
    
The reversed game would be played by starting at $1$ and adding either $1, 2, $ or $3$ (never going above $n$), with the first person to make the number $n$ winning. 

This is transparently the same as Nim, with the same winning strategy for Player $2$. Thus, this game does not demonstrate any new behavior when reversed: The winning structure is still the same as Player $2$ always wins with the same strategy.
\end{example}
In general, reversed games do not exhibit the same winning structure or winning strategies. As we will see, reversed Chomp has a similar winning structure but different winning strategy, and the reversed Zeckendorf game has an entirely different winning structure. 
%=======================================
%=========  reversed chomp  ============
%=======================================
\subsection{Reversed Chomp}
As a useful example of a reversed game, we solve the reversed Chomp game. First, we introduce the game of \emph{Chomp}. Chomp is a normal, impartial game played on a rectangular board with $N$ rows and $M$ columns. Players take turns choosing a square on the board and ``eating" it, removing that square and all the squares both above and to the right of it. The player forced to eat the ``poisoned square" in the bottom left corner loses the game. We consider the game played with two players and ignore the trivial game $N=M=1$.

We can consider this game as played on a directed graph with a starting node representing the full board and an ending node representing the poisoned square. It can be shown by a ``strategy-stealing argument'' that for all $N,M$, Chomp is a win for Player 1. However, constructing an explicit winning strategy proves difficult.

Consider \emph{reversed Chomp}, which starts from the poisoned square and ends with the full board. Without loss of generality, assume $N\geq M$. When $M=1$, Player 1 wins in one move. 

We provide a constructive proof that $M>1$ is always a win for Player 2 by giving the following winning strategy.

\begin{enumerate}[label = (\roman*)]
\item \textbf{Case $M=2$}: Player 1 has three options for their first move: complete the bottom row, complete the leftmost column, or partially complete the leftmost column. The first two choices allow Player 2 to win in one move. 
    
Suppose Player 1 partially completed the leftmost column to a height of $1<h<N$ squares. Then Player 2 can respond by filling in the other column to a height of $h-1$. The game has now been reduced to another game with 2 columns but with fewer rows than we started with (see figure \hyperref[fig:case_M_equals_2]{1}).

Player 1 is presented with the same three choices. Player 2 can continue responding in the above manner until we reach the game with 2 rows and 2 columns. From this position, regardless of Player 1's choice, Player 2 wins in one move.
\item \textbf{Case $M>2$}: Player 1 has four options for their first move: complete the bottom row, complete the leftmost column, partially complete the bottom row, or partially complete the leftmost column. The first two choices allow Player 2 to win in one move.

If Player 1 partially completes the bottom row to a length of $1<l<M$ squares, then Player 2 can respond by filling in the leftmost $l-1$ columns (see figure \hyperref[fig:case_M_greater_than_2]{2}). This reduces to a game with the same number of rows and fewer columns than we started with. Similarly, if Player 1 partially completes the leftmost column to a height of $1<h<N$ squares, then Player 2 can respond by filling in the bottom $h-1$ rows which reduces to a game with the same number of columns but fewer rows (see Figure \hyperref[fig:case_M_greater_than_2_image2]{3}).

Either way, Player 1 is presented with the same four choices. Player 2 can continue responding in this manner until the game is reduced to a state with either only 2 rows or 2 columns. By the previous case, this is a win for Player 2.
\end{enumerate}

%%%%%%%%%%%%%%%%%%%%%%%%%%%%%%%%%%%%%%%%%%%

\setcounter{figure}{0}

\begin{figure}
    \label{fig:case_M_equals_2}
\begin{tikzpicture}
    % Draw grid lines
    \draw[step=1cm, gray, very thin] (0,0) grid (2,5);
    
    % Draw poisoned square
    \fill[red] (0,0) rectangle (1,1);
    
    % Label for poisoned square
    \node[white] at (0.5,0.5) {X};
    
    % Black squares
    \fill[gray] (0,1) rectangle (1,3);
    \fill[black] (1,0) rectangle (2,2);
\end{tikzpicture}
\hspace{1cm}
$\equiv$
\hspace{1cm}
\begin{tikzpicture}
    % Draw grid lines
    \draw[step=1cm, gray, very thin] (0,0) grid (2,3);
    
    % Draw poisoned square
    \fill[red] (0,0) rectangle (1,1);
    
    % Label for poisoned square
    \node[white] at (0.5,0.5) {X};
\end{tikzpicture}
\caption{$M=2$. Gray indicates move by Player 1, black indicates move by Player 2}
\end{figure}
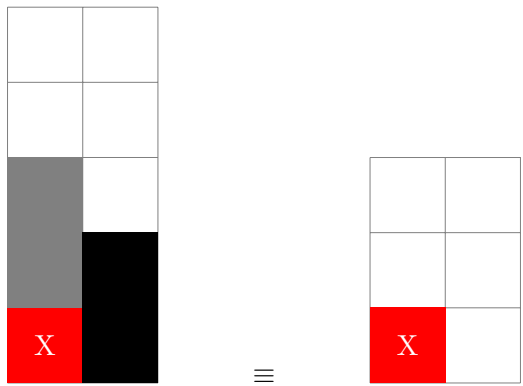

%%%%%%%%%%%%%%%%%%%%%%%%%%%%%%%%%%%%%%%%%%%%%%%%%%%%%%%%%%%%%%%

\begin{figure}
    \label{fig:case_M_greater_than_2}
\begin{tikzpicture}
    % Draw grid lines
    \draw[step=1cm, gray, very thin] (0,0) grid (5,5);
    
    % Draw poisoned square
    \fill[red] (0,0) rectangle (1,1);
    
    % Label for poisoned square
    \node[white] at (0.5,0.5) {X};
    
    % Black squares
    \fill[gray] (1,0) rectangle (3,1);
    \fill[black] (0,1) rectangle (2,5);
\end{tikzpicture}
\hspace{1cm}
$\equiv$
\hspace{1cm}
\begin{tikzpicture}
    % Draw grid lines
    \draw[step=1cm, gray, very thin] (0,0) grid (3,5);
    
    % Draw poisoned square
    \fill[red] (0,0) rectangle (1,1);
    
    % Label for poisoned square
    \node[white] at (0.5,0.5) {X};
\end{tikzpicture}
\caption{$M>2$, Player 1 partially completes bottom row.}
\end{figure}
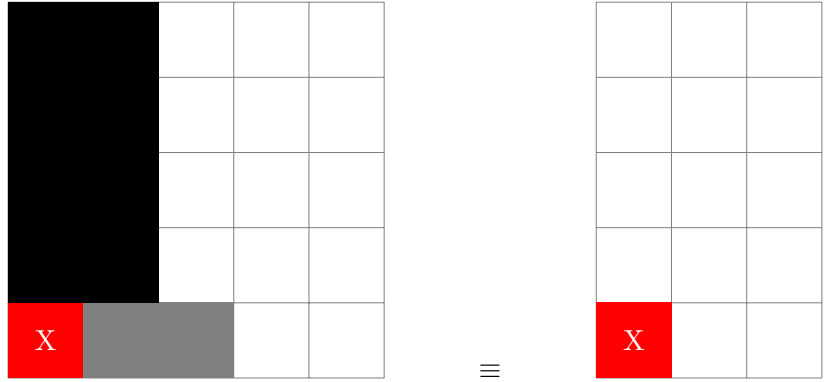

%%%%%%%%%%%%%%%%%%%%%%%%%%%%%%%%%%%%

\begin{figure}
    \label{fig:case_M_greater_than_2_image2}
\begin{tikzpicture}
    % Draw grid lines
    \draw[step=1cm, gray, very thin] (0,0) grid (5,5);
    
    % Draw poisoned square
    \fill[red] (0,0) rectangle (1,1);
    
    % Label for poisoned square
    \node[white] at (0.5,0.5) {X};
    
    % Black squares
    \fill[gray] (0,1) rectangle (1,4);
    \fill[black] (1,0) rectangle (5,3);
\end{tikzpicture}
\hspace{1cm}
$\equiv$
\hspace{1cm}
\begin{tikzpicture}
    % Draw grid lines
    \draw[step=1cm, gray, very thin] (0,0) grid (5,2);
    
    % Draw poisoned square
    \fill[red] (0,0) rectangle (1,1);
    
    % Label for poisoned square
    \node[white] at (0.5,0.5) {X};
\end{tikzpicture}
\caption{$M>2$, Player 1 partially completes leftmost column.}
\end{figure}
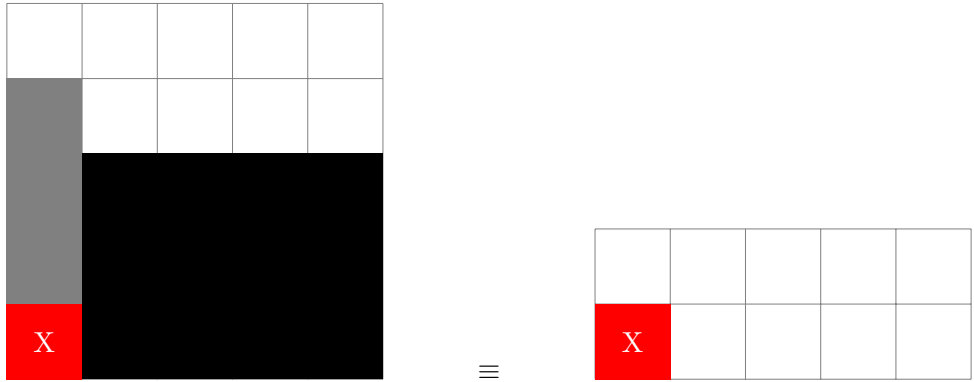

\begin{remark}
    Along with the winning strategy, the winning structure for reversed Chomp is also different than normal Chomp. Player 1 wins reversed Chomp for all $N$ when $M=1$. This is an infinite number of games but not a positive proportion. In forward Chomp, Player 2 only wins a finite number of games, namely just the trivial game.
\end{remark}

\section{The Reversed Zeckendorf Game}\label{sec:reversed_zeck_game}

Following our definition of the reversed game in Definition \ref{def:reverse}, we give an explicit construction for the \emph{reversed Zeckendorf game.} 

The terminology is the same as in the forwards game (with the starting and ending positions flipped, so we always start at the Zeckendorf decomposition of $n$) except we rename each move.
\begin{enumerate}[label = (\roman*)]
    \item \textit{Split}: We remove one chip from the $(i+1)$\textsuperscript{th} bin and place one chip each in the $i$\textsuperscript{th} and $(i-1)$\textsuperscript{th} bins.
    \begin{align}
        F_{i+1}&\ \mapsto\  F_{i-1}+ F_{i}, \\
        F_{2} &\ \mapsto\ 2F_1. \nonumber
    \end{align}
    \item \textit{Combine}: For $i> 2$, we have
    \begin{align}
        F_{i-2} + F_{i+1} &\ \mapsto\  2F_i, \\
        F_{3} + F_{1} &\ \mapsto\ 2F_2. \nonumber
    \end{align}
\end{enumerate}

\PlayerOneWinsInfinite*
We present two proofs, one constructive and the other nonconstructive.

\begin{proof}(Nonconstructive)
By way of contradiction, suppose Player $2$ has a winning strategy. If Player 1 chose to combine for their first move, Player 2 has a forced win starting at the state $2F_{i}$. There is only one move in this position, which means Player 2 has a forced win with Player 1 starting at the state $F_i+F_{i-1} + F_{i-2}$. 

However, Player 1 could choose to instead split the $(i+1)$\textsuperscript{th} bin, which makes it Player 2's turn at the state $F_i+F_{i-1} + F_{i-2}$. Now, Player 1 can steal the strategy from Player 2 to have a forced win in the starting state, a contradiction.
\end{proof}

\begin{remark}
The constructive proof relies on Lemma \ref{lem:evenwin}. As a corollary of this lemma, we have a constructive proof for why the game state $2F_i$ is a win for whoever goes second. Before that, Player 1 starts with combining to bring the game into this situation.
\end{remark}

The next natural question is whether or not the same can be said of Player 2.

\PlayerTwoWinsToo*

To investigate this conjecture, we directly computed which player had a forced winning strategy for $n\leq 129$ (see Appendix \ref{sec:table}). The code for computing which player has a forced winning strategy is listed in Appendix \ref{sec:code} and has a algorithmic complexity of $O(\exp(\sqrt{n})$. For $n = 129$, the program took close to $2$ hours. 

\begin{figure}
    \hspace*{5em}
    \includegraphics[scale=0.9]{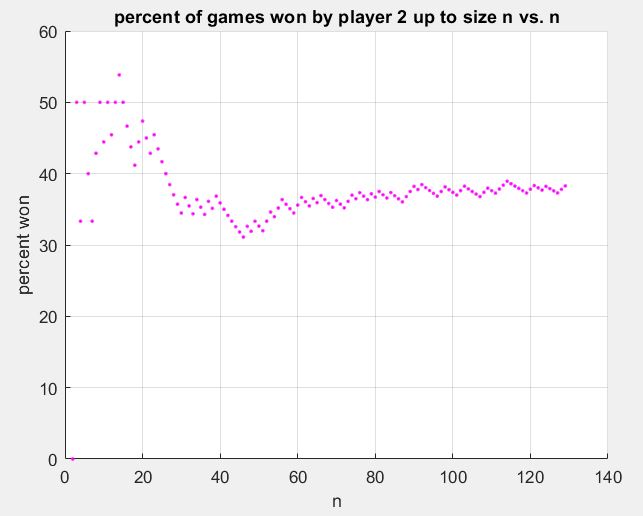}
    \label{fig:percent_player_two_wins_reverse_zeck}
    \caption{Percent of games $\leq n$ won by Player 2}
\end{figure}

Figure \hyperref[fig:percent_player_two_wins_reverse_zeck]{4} shows the proportion Player 2 wins. We plot $n$ versus the percent of games won by Player 2 in games $2$ through $n$.

For the first $n\leq 129$, the average number of Player 1 wins is $80/129 \approx .620$. The average number of wins seems to stabilize computationally. Combining these two together, we are led to the following conjecture.

\WinningPercentage*

This is a natural conjecture due to the connection between the Fibonacci numbers and the golden ratio. Work on determining the winner is provided in \textsection \ref{sec:startpos}. 

\subsection{Other Facts}\label{sec:otherfacts}
We gather various results  about the reversed Zeckendorf game which follow from already-known facts about the forwards game. Since we have merely reversed the arrows in the forward game tree to produce the reversed game, many properties about the forwards game extend immediately as corollaries.
One property that extends is the lengths of games. From \cite{original_paper} for the lower bound and \cite{TowardsGauss} for the upper bound, we have the following.

\gamelengths*

Theorems 1.8 and 1.9 from the paper \cite{TowardsGauss} also transfer. 

Finally, we have results on randomly played games with one subtlety. The paper discusses two different probability measures put on the space of all possible games with a fixed $n$. One of them simply assigns the same probability to each game (denoted $\mu_n$), as well as assigning each game according to the probability of playing it if each player picked a random move each turn, uniformly from the possible moves (denoted $\mathbb{P}_n$).

Now, the first measure $\mu$ assigns the same probability to each game as its corresponding flipped game in the reversed version. As such, statements about this probability measure easily translate as a corollary. However, the measure $\mathbb{P}_n$ requires more work than deducing the result for the reversed game from the original game. In this case, modifications and checks are required throughout Section $4$ of their paper until Lemma $4.5$. In the end, the main result still holds.

\randomgames*

This theorem works for a more general setting than that of the two player game. The game can be modified to $Z$ many players, wherein the players take turns in sequence and the last player who moves wins. The theorem can then be interpreted as saying that in the $Z$ player game, where each player chooses their move uniformly at random from their possible choices, any given player will win approximately $1/Z$ percent of the time.

For the specific case of two players and an even $n$, we have the following.
\begin{theorem} For the Reversed Zeckendorf Game player on an even $n$,
\begin{equation}
\mathbb{P}_n(\text{game length equals }1 \text{ mod }2)\ =\ \frac{1}{2}.
\end{equation}
\end{theorem}
\begin{proof}
Equivalently, we can show that each player, playing randomly on each turn, each has a $1/2$ probability of winning.
    
At some point, every game will have the $h_3 = 1$ and $h_i = 0$ for all $i > 3$. Furthermore, we can continue to the move that gets rid of the three. We have
\begin{equation}
n\ =\ h_1 + 2h_2+3. 
\end{equation}
Since $n$ is even, this means $h_1 > 0$, so at the final turn $h_3 = 1$, the player has 2 choices: (i) combine a one and a three or (ii) split the three.

Once there are no threes (i.e., $h_3 = 0$), the game only has one option, namely splitting the two until there are no moves left. The winner is determined by the parity of the second bin. However, (i) will yield the height of bin 2 equals $h_2 + 2$ and (ii) will yield the height of bin 2 equals $h_2 + 1$. 

There is an equal chance the player will lose or win, and nothing that came before or after matters. We conclude the probability either player wins is $1/2$, as desired.
\end{proof}

%==============================================
%============ OTHER STARTING POSITIONS ============
%==============================================

\section{Other Starting Positions}\label{sec:startpos}

A natural modification of the reversed Zeckendorf game is to alter the starting position. Instead of starting at the Zeckendorf decomposition, we instead specify a starting position and start the game from there. In fact, such problems are often tractable as the starting and ending positions are well-understood. Contrast this with applying the same idea to the forwards game: In this case, we would still be headed towards the Zeckendorf decomposition, which can vary as $n$ changes.
\begin{lemma}\label{lem:evenwin}
If the starting position of the reversed Zeckendorf game has all bins of even height, then Player 2 has a winning strategy.
\end{lemma}
\vspace{-5mm}
\begin{proof}
The proof uses a ``copycat'' strategy that is, if Player 1 does a move, Player 2 plays the exact same move. We must show the following:
\begin{enumerate}
\item if Player 1 makes a move, then Player 2 can also make that move, and
\item after Player 2 moves, all the bins are of even height.
\end{enumerate}
If these are both the case, then Player 1 will have to move on a position where there are only ones and twos. At this point, there will be an even number of twos, so the game from here on is deterministic; both players will split twos until there are no twos left. Since the parity of the twos is even, this means it will be Player 1's turn when there are no twos left, so they lose.

Consider Player 1 making a split move so that \begin{align*}
h_i &\ \mapsto\ h_i - 1,  \\
h_{i-1} &\ \mapsto\ h_{i-1} + 1, \text{ and}\\
h_{i-2} &\ \mapsto\ h_{i-2} + 1.
\end{align*}
Since $h_i$ is even, then $h_i - 1$ is odd. Thus, $h_i - 1 \geq 1 > 0$. This means Player 2 can legally make the same splitting move with the net effect of   
    \begin{align*}
        h_i &\ \mapsto\ h_i - 2, \\
        h_{i-1} &\ \mapsto\ h_{i-1} + 2, \text{ and}\\
        h_{i-2} &\ \mapsto\ h_{i-2} + 2,
    \end{align*}
    with all the other bins unchanged, so all the heights are still even. We can do a similar argument for the combine move which proves (1) and (2).
\end{proof}
We also solve the game when the starting position consists only of ones, twos, and threes. We denote the starting position with $a$ ones, $b$ twos, and $c$ threes as the ordered triple $(a, b, c)$ from this point on. We will also often use the notation that an $O$ or $E$ in the place of the element of the tuple means that the corresponding variable is odd or even respectively (i.e., $(E, E, E)$ means that $a, b$ and $c$ are all even). Further, $O'$ and $E'$ mean that the corresponding variable is odd or even, but independent from the previous $O$ and $E$.

\OneTwoThree*
\begin{proof}
The first case is a corollary of Lemma \ref{lem:evenwin}, since all the bins have even height. The next three cases follow from the first case as Player 1 simply splits the 3, splits the 2, or combines the 1 and the 3 respectively to bring the game state into all three bins being even. 

The next few cases are a bit more work. We show the (Odd, Odd, Even) case, where $a>c$, and put the rest in Appendix \ref{sec:proofthm123}. 

For this case, we assume $a$ is odd and $b$ and $c$ are even with $a>c$. We show Player 2 wins by force.

First, assume Player 1 splits one of the threes. This gives the game state
\begin{equation}
    (a+1, b+1, c-1).
\end{equation}
Since $c$ is even, we have $c-1 \geq 1$. This allows Player 2 to split another three, yielding
\begin{equation}
    (a+2,b+2,c-2).
\end{equation}
Note that we have $a+2 > c-2$. Moreover, $a-2$ is odd while $b+2$ and $c-2$ remain even. By induction, Player 2 wins.

Next, assume Player 1 splits one of the twos. This yields the game state
\begin{equation}
    (a+2, b-1, c).
\end{equation}
Since $b$ is even, we have $b-1 \geq 1$. This allows Player 2 to split another two, yielding
\begin{equation}
    (a+4, b-2,c).
\end{equation}
Note we have $a+4 > c$. Moreover, $a+4$ is odd while $b-2$ and $c$ remain even. By induction, Player 2 wins.

Finally, assume Player 1 combines a one and a three. This yields the game state
\begin{equation}
    (a-1, b+2, c-1).
\end{equation}
Since $c$ is even, we have $c-1 \geq 1$. Moreover, since $a > c$ by assumption, we have $a > 1$. In particular, Player 2 is able to perform another combine, leading to
\begin{equation}
    (a-2, b+4, c-2).
\end{equation}
Note $a-2 > c-2$. Moreover, all parities are preserved. By induction, Player 2 wins.
\end{proof}

\begin{remark}
We summarize two key points in the above proof.
\begin{enumerate}[label = (\roman*)]
\item The method in the above Theorem is constructive, so following it carefully provides a winning strategy.
\item The last two cases appear to not depend on whether $a>c$ or $a<c$, but the strategy to win changes in those cases. 
\end{enumerate}
\end{remark}

\section{The Build Up Game}\label{sec:buildup}
We move on to the \emph{Build Up 1-2-3 game}.
In this game, two players, say Player 1 and Player 2, take turns putting down a 1, 2 or 3 until their sum equals exactly $n$, generating an ordered triple $(a, b, c)$. After reaching exactly $n$, the players start playing the reverse Zeckendorf game on this ordered triple and start with the player who did not move last. Whoever wins this reversed Zeckendorf game wins the whole game. We solved this game.
\begin{restatable}{theorem}{BuildUpOneTwoThree}\label{thm:buildup123}
For $n \neq 4$,
\begin{align}
        n = 4 &\implies \text{Player 1 wins}, \\
        n \text{ odd} &\implies \text{Player 1 wins}, \\
        n \neq  4 \text{ even }  &\implies \text{Player 2 wins}.
    \end{align}
\end{restatable}
\begin{proof}
Using Theorem \ref{thm:123}, we know who wins in each ordered triple. We provide a constructive proof in which we split the game into cases mod 4. We consider residue classes modulo 4 since any player may play 4 minus what was played previously and thus preserve parity modulo 4. This will be referred to as \emph{nimming down}. An extra property is that the number of ones and threes put down is the same, and the number of twos is even.

For $n \equiv 0 \: (4)$ and $n \geq 8$, Player 2 should nim down until there is 8 left to play. From here, there are two cases.

\begin{enumerate}[label = (\arabic*)]
    \item Player 1 is left to play on $(E, E', E)$:
    \begin{enumerate}[label = (\roman*)]
        \item If Player 1 puts down a 1 or a 3, then Player 2 puts down the other, leaving $(E+1, E', E+1)$.
    
        \begin{enumerate}[label = (\alph*)]
            \item If Player 1 puts down a 1 or a 3, then Player 2 puts down the other, leaving player 1 to start on $(E+2,E',E+2)$, so Player 2 wins.
            \item If Player 1 puts down a 2, then Player 2 puts down a 1, forcing Player 1 to put down a 1. Player 2 then starts on $(E+3,E'+1,E+1)$, which is a win for Player 2.
        \end{enumerate}
        
        \item If Player 1 puts a 2 down, Player 2 puts down a 1, leaving Player 1 to play on $(E+1,E'+1,E)$, with 5 total left to play.
    
\begin{enumerate}[label = (\alph*)]
\item If Player 1 puts down a 1 or a 3, Player 2 puts down the other, forcing Player 1 to put down a 1. This leaves Player 2 to start on $(E+3,E'+1,E+1)$, which is a win for Player 2.
\item If Player 1 puts down a 2, then Player 2 puts down a 2, and Player 1 is forced to put down a 1. Player 2 then starts on $(E+2,E'+3,E)$, which is a win for Player 2. 
        \end{enumerate}
    \end{enumerate}
    \item Player 1 is left to play on $(O,E,O)$:
    \begin{enumerate}[label = (\roman*)]
        \item If Player 1 puts down a 1, then Player 2 puts down a 2, so Player 1 plays on $(O+1,E+1,O)$, with a total of 5 left to play. 
        \begin{enumerate}[label = (\alph*)]
            \item If Player 1 puts down a 1 or 3, then Player 2 puts down the other, forcing Player 1 to put down a 1. This leaves Player 2 to start on $(O+3,E+1,O+1)$, which is a win for Player 2. 
            \item If Player 1 puts down a 2, then Player 2 puts down a 2, forcing Player 1 to put down a 1. Player 2 starts on $(O+2,E+3,O)$, which is a win for Player 2.
        \end{enumerate}
        \item  If Player 1 puts down a 2 or 3, then Player 2 puts down the other, so Player 1 plays on $(O,E+1,O+1)$ with 3 total left.  
        \begin{enumerate}[label = (\alph*)]
            \item If Player 1 puts down a 3, then Player 2 starts the game on $(O,E+1,O+2)$, which is a win for Player 2.
            \item If Player 1 puts down a 1 or 2, then Player 2 puts down the other. Player 1 then starts the game on $(O+1,E+2,O+1)$, which is a win for Player 2. 
        \end{enumerate}
    \end{enumerate}
\end{enumerate}
The other three cases are very similar. These, as well as the winning strategy for small values of $n$, have been included in the appendix. 
\end{proof}

\section{Further Directions and Conclusions}
There are numerous directions for future work.
\begin{enumerate}[label = (\roman*)]
\item While the authors believe Conjecture \ref{conj:strong} is very much out of reach, Conjecture \ref{conj:infinitewin} can be proven by finding a family of Zeckendorf decompositions where Player $2$ can force a win. 
\item More generally, one might study other reversed games to see if any other games reveal such a complex winning structure upon reversal. 
\item Future work could also try to solve the Reversed game for a larger family of starting positions, perhaps ones closer to the Zeckendorf decomposition.

\item Other modifications of the game can be considered, including a ``stagnant one" variation, where all chips with value one are removed from the game (or equivalently, any move requiring a one cannot be played).
\end{enumerate}

We conclude with a final, bold conjecture.
\begin{conjecture}
    The winning strategy for the forward Zeckendorf games is dependent upon who wins the reverse Zeckendorf game.
\end{conjecture}
This is a potential explanation of why a constructive proof of why Player $2$ always wins the forwards game has remained elusive. 

%%%%%%%%%%%%%%%%%%%%%%%%%%%%%%%%%%%%%%%%%%%%%%%%%%%%%%%%%%%%%%%%%%%%%%%%%%%%%%%%%%%%%%%%%%%%%%%%%%%%%%%%%%%%%%%%%%%%%%%%
%%%%%%%%%%%%%%%%%%%%%%%%%%%%%%%%%%%%%%%%%%%%%%%%%%%%%%%%%%%%%%%%%%%%%%%%%%%%%%%%%%%%%%%%%%%%%%%%%%%%%%%%%%%%%%%%%%%%%%%%
%%%%%%%%%%%%%%%%%%%%%%%%%%%%%%%%%%%%%%%%%%%%%%%%%%%%%%%%%%%%%%%%%%%%%%%%%%%%%%%%%%%%%%%%%%%%%%%%%%%%%%%%%%%%%%%%%%%%%%%%
\newpage

\appendix

\section{Additional Details}

\subsection{Table of Wins for Reverse Zeckendorf Game}\label{sec:table}
The following table lists which player wins for the starting number $n$ as well as the number of vertices and edges of the associated directed graph. A grayed square means Player 2 wins while a white one means Player 1 wins.

\begin{center}
\renewcommand{\arraystretch}{1.2}{
\begin{tabular}{|P{0.5cm}|P{1.5cm}|P{1.5cm}|P{1.5cm}| }
\hline
$n$ & Result & Edges & Vertices \\
\hline
2 & & 1 & 2 \\
\hline
3 & \cellcolor{lightgray} & 2 & 3 \\
\hline
4 & & 4 & 4 \\
\hline
5 & \cellcolor{lightgray} & 7 & 6 \\
\hline
6 & & 11 & 8 \\
\hline
7 & & 16 & 10 \\
\hline
8 & \cellcolor{lightgray} & 24 & 14 \\
\hline
9 & \cellcolor{lightgray} & 32 & 17 \\
\hline
10 & & 45 & 22 \\
\hline
11 & \cellcolor{lightgray} & 59 & 27 \\
\hline
12 & & 77 & 33 \\
\hline
13 & \cellcolor{lightgray} & 100 & 41 \\
\hline
14 & \cellcolor{lightgray} & 126 & 49 \\
\hline
15 & & 158 & 59 \\
\hline
16 & & 198 & 71 \\
\hline
17 & & 241 & 83 \\
\hline
18 & & 297 & 99 \\
\hline
19 & \cellcolor{lightgray} & 358 & 115 \\
\hline
20 & \cellcolor{lightgray} & 430 & 134 \\
\hline
21 & & 516 & 157 \\
\hline
22 & & 610 & 180 \\
\hline
23 & \cellcolor{lightgray} & 722 & 208 \\
\hline
24 & & 849 & 239 \\
\hline
25 & & 990 & 272 \\
\hline
26 & & 1158 & 312 \\
\hline
27 & & 1339 & 353 \\
\hline
28 & & 1548 & 400 \\
\hline
29 & & 1785 & 453 \\
\hline
30 & & 2043 & 509 \\
\hline
31 & \cellcolor{lightgray} & 2341 & 573 \\
\hline
32 & & 2667 & 642 \\
\hline
33 & & 3028 & 717 \\
\hline
\end{tabular}}
\renewcommand{\arraystretch}{1.2}{
\begin{tabular}{|P{0.5cm}|P{1.5cm}|P{1.5cm}|P{1.5cm}| }
\hline
$n$ & Result & Edges & Vertices \\
\hline
34 & \cellcolor{lightgray} & 3440 & 803 \\
\hline
35 & & 3881 & 892 \\
\hline
36 & & 4381 & 993 \\
\hline
37 & \cellcolor{lightgray} & 4930 & 1102 \\
\hline
38 & & 5528 & 1219 \\
\hline
39 & \cellcolor{lightgray} & 6199 & 1350 \\
\hline
40 & & 6924 & 1489 \\
\hline
41 & & 7721 & 1640 \\
\hline
42 & & 8606 & 1808 \\
\hline
43 & & 9552 & 1983 \\
\hline
44 & & 10606 & 2178 \\
\hline
45 & & 11743 & 2386 \\
\hline
46 & & 12979 & 2609 \\
\hline
47 & \cellcolor{lightgray} & 14339 & 2854 \\
\hline
48 & & 15796 & 3113 \\
\hline
49 & \cellcolor{lightgray} & 17387 & 3393 \\
\hline
50 & & 19119 & 3697 \\
\hline
51 & & 20970 & 4017 \\
\hline
52 & \cellcolor{lightgray} & 23001 & 4367 \\
\hline
53 & \cellcolor{lightgray} & 25171 & 4737 \\
\hline
54 & & 27517 & 5134 \\
\hline
55 & \cellcolor{lightgray} & 30064 & 5564 \\
\hline
56 & \cellcolor{lightgray} & 32777 & 6016 \\
\hline
57 & & 35719 & 6504 \\
\hline
58 & & 38879 & 7025 \\
\hline
59 & & 42252 & 7575 \\
\hline
60 & \cellcolor{lightgray} & 45909 & 8171 \\
\hline
61 & \cellcolor{lightgray} & 49794 & 8797 \\
\hline
62 & & 53971 & 9466 \\
\hline
63 & & 58458 & 10183 \\
\hline
64 & \cellcolor{lightgray} & 63222 & 10936 \\
\hline
65 & & 68351 & 11744 \\
\hline
\end{tabular}}
\end{center}

\newpage
\noindent

\begin{center}
\renewcommand{\arraystretch}{1.2}{
\begin{tabular}{|P{0.5cm}|P{1.5cm}|P{1.5cm}|P{1.5cm}| }
\hline
$n$ & Result & Edges & Vertices \\
\hline
66 & \cellcolor{lightgray} & 73811 & 12599 \\
\hline
67 & & 79627 & 13502 \\
\hline
68 & & 85874 & 14471 \\
\hline
69 & & 92480 & 15486 \\
\hline
70 & \cellcolor{lightgray} & 99552 & 16568 \\
\hline
71 & & 107083 & 17715 \\
\hline
72 & & 115060 & 18921 \\
\hline
73 & \cellcolor{lightgray} & 123593 & 20207 \\
\hline
74 & \cellcolor{lightgray} & 132622 & 21559 \\
\hline
75 & & 142212 & 22987 \\
\hline
76 & \cellcolor{lightgray} & 152437 & 24506 \\
\hline
77 & & 163215 & 26094 \\
\hline
78 & & 174701 & 27782 \\
\hline
79 & \cellcolor{lightgray} & 186852 & 29558 \\
\hline
80 & & 199697 & 31425 \\
\hline
81 & \cellcolor{lightgray} & 213356 & 33405 \\
\hline
82 & & 227747 & 35478 \\
\hline
83 & & 242990 & 37664 \\
\hline
84 & \cellcolor{lightgray} & 259136 & 39973 \\
\hline
85 & & 276121 & 42386 \\
\hline
86 & & 294140 & 44939 \\
\hline
87 & & 313108 & 47613 \\
\hline
88 & \cellcolor{lightgray} & 333117 & 50421 \\
\hline
89 & \cellcolor{lightgray} & 354284 & 53384 \\
\hline
90 & \cellcolor{lightgray} & 376512 & 56478 \\
\hline
91 & & 399993 & 59735 \\
\hline
92 & \cellcolor{lightgray} & 424730 & 63154 \\
\hline
93 & & 450710 & 66727 \\
\hline
94 & & 478155 & 70492 \\
\hline
95 & & 506942 & 74422 \\
\hline
96 & & 537240 & 78543 \\
\hline
97 & \cellcolor{lightgray} & 569148 & 82871 \\
\hline
\end{tabular}}
\renewcommand{\arraystretch}{1.2}{
\begin{tabular}{|P{0.5cm}|P{1.5cm}|P{1.5cm}|P{1.5cm}| }
\hline
$n$ & Result & Edges & Vertices \\
\hline
98 & \cellcolor{lightgray} & 602577 & 87383 \\
\hline
99 & & 637787 & 92122 \\
\hline
100 & & 674725 & 97075 \\
\hline
101 & & 713458 & 102247 \\
\hline
102 & \cellcolor{lightgray} & 754214 & 107677 \\
\hline
103 & \cellcolor{lightgray} & 796847 & 113331 \\
\hline
104 & & 841620 & 119251 \\
\hline
105 & & 888582 & 125442 \\
\hline
106 & & 937699 & 131890 \\
\hline
107 & & 989279 & 138644 \\
\hline
108 & \cellcolor{lightgray} & 1043215 & 145681 \\
\hline
109 & \cellcolor{lightgray} & 1099679 & 153022 \\
\hline
110 & & 1158887 & 160703 \\
\hline
111 & & 1220686 & 168686 \\
\hline
112 & \cellcolor{lightgray} & 1285449 & 177031 \\
\hline
113 & \cellcolor{lightgray} & 1353148 & 185727 \\
\hline
114 & \cellcolor{lightgray} & 1423848 & 194777 \\
\hline
115 & & 1497885 & 204232 \\
\hline
116 & & 1575109 & 214059 \\
\hline
117 & & 1655815 & 224299 \\
\hline
118 & & 1740180 & 234978 \\
\hline
119 & & 1828092 & 246065 \\
\hline
120 & \cellcolor{lightgray} & 1920015 & 257632 \\
\hline
121 & \cellcolor{lightgray} & 2015842 & 269652 \\
\hline
122 & & 2115775 & 282150 \\
\hline
123 & & 2220148 & 295175 \\
\hline
124 & \cellcolor{lightgray} & 2328793 & 308687 \\
\hline
125 & & 2442154 & 322751 \\
\hline
126 & & 2560320 & 337374 \\
\hline
127 & & 2683290 & 352543 \\
\hline
128 & \cellcolor{lightgray} & 2811589 & 368337 \\
\hline
129 & \cellcolor{lightgray} & 2945040 & 384715 \\
\hline
\end{tabular}}
\end{center}

\subsection{Proof of Theorem \ref{thm:123}}\label{sec:proofthm123}

\OneTwoThree*

\begin{proof}
For the $(O,E,E)$ case, where $a<c$, we show Player 1 wins by force.

Player 1 wins by performing a combine with a one and a three, yielding the game state
\begin{equation}
(a-1, b+2, c-1).
\end{equation}
Note $a-1 < c-1$. This also yields the $(E,E,O)$ case with Player 2 on move. By work above, the player on move loses which means Player 1 wins.

For the $(E,E,O)$ and $a>c$ case, we show Player 1 wins by force.

Player 1 wins by combining a one and a three. This gives the game state
\begin{equation}
(a-1, b+2, c-1).
\end{equation}
Then, $a-1$ is odd, $b+2$ is even, and $c-1$ is even. Moreover, we have $a-1 > c-1$. Therefore, by the $(O,E,E)$ proof, the player on move loses. Since Player 2 is on move, then Player 1 wins.

For the $(E,E,O)$ and $a<c$ case, we show Player 2 wins by force.

First, assume Player 1 combines a one and a three. This gives the game state
\begin{equation}
(a-1, b+2, c-1).
\end{equation}
Since Player 1 was able to combine a one and a three, then $a\geq 2$ which implies $c > 2$. We have $a-1, c-1 \geq 1$, and so Player 2 may perform another combine to get
\begin{equation}
(a-2, b+4, c-2).
\end{equation}
The point is that $a-2 < c-2$ and all parities are preserved. Eventually, Player 1 will run out of ones and threes to combine and is forced to either split a three or a two. Both of these cases are covered next, and are wins for Player 2. (Informally, the main idea in this case is to force Player 1 to run out combines, upon which they will be forced to perform a split. We will show that all splits are losing for Player 1, and so Player 2 forces a win by running Player 1 out of combines.)

Assume Player 1 splits one of the twos. This gives the game state
\begin{equation}
(a+2,b-1,c).
\end{equation}
Now, if $a+3 > c-1$, Player 2 wins by splitting one of the threes. (Since $c$ is odd, a three exists.) This yields the game state
\begin{equation}
    (a+3, b, c-1).
\end{equation}
At this point, $a+3$ is odd, $b$ is even, and $c-1$ is even. Moreover, we have $a+3 > c-1$. By work above, the player on move loses. Since Player 1 is on move, Player 2 wins.

On the other hand, if $a+3 < c-1$ (i.e., $a+4 < c$), Player 2 should copy Player 1 and split another two. Note that a two exists; since $b$ is even and Player 1 was able to split a two, we must have $b \geq 2$. This forces $b-1 \geq 1$, and so a two is available. This yields the game state
\begin{equation}
(a+4, b-2, c).
\end{equation}
The upshot is the number of ones is increasing, and so Player 2 will eventually be able to reduce to the situation analyzed in the $a+3 > c-1$ case. The parities are preserved and $a+4< c$. Therefore, Player 2 wins.

Finally, assume Player 1 splits one of the threes. This yields the game state
\begin{equation}
(a+1, b+1, c-1).
\end{equation}
If $c = 1$, then $c-1=0$. Since $a<c$, we get $a = 0$. The point is that we have the game state $(1, b+1, 0)$ with Player 2 on move. The only option for both players at this point is to keep splitting twos. Since $b+1$ is odd, Player 1 will run out of splits first, and so Player 2 wins.

Now, assume $c > 1$. If $a+3>c-1$, Player 2 wins by splitting a two, yielding the position
\begin{equation}
(a+3, b, c-1).
\end{equation}
Note $a+3$ is odd, $b$ is even, and $c$ is even. Since $a+3 > c-1$, then earlier work implies the player on move loses. As Player 1 is on move, Player 2 wins.

Finally, if $a+3 <c-1$, Player 2 wins by splitting a three, yielding the position
\begin{equation}
    (a+2, b+2, c-2).
\end{equation}
All parities are preserved, and $a+2 < c-2$ since $a+3 < c-1$. By induction, Player 2 wins.

In all cases, we have shown Player 2 has a forced win.

For the $(E,O,O)$ case, we show Player 1 wins by force.

The winning strategy depends on the number of ones versus the number of threes. First, assume $a + 2 > c$; that is, there is a large number of ones. Player 1 wins by splitting a three, yielding
\begin{equation}
(a+1, b+1, c-1).
\end{equation}
Note $a+1$ is now odd, while $b+1$ and $c-1$ are even. Furthermore, we have $a+1 > c-1$ since $a+2 > c$. By earlier work, the player on move loses. Since Player 2 is on move, then Player 1 wins.

On the other hand, we assume $a + 2 < c$; that is, there is a small number of ones. Player 1 wins by splitting a two, yielding the game state
\begin{equation}
(a+2, b-1, c).
\end{equation}
Note $a+2$ is even, $b-1$ is even, and $c$ is odd. Moreover, we have $a+2 < c$. By earlier work, the player on move loses. Since Player 2 is on move, then Player 1 wins.

In either case, we have exhibited a winning strategy for Player 1, and so Player 1 wins by force.

For the $(O,O,E)$ case, we assume $a$ and $b$ are odd and $c$ is even. We show Player 1 wins by force.

As before, the winning strategy depends on the number of ones versus the number of threes. First, assume $a+2 < c$. Player 1 wins by splitting a three, yielding
\begin{equation}
    (a+1, b+1, c-1).
\end{equation}
Note $a+1$ and $b+1$ are now even, while $c-1$ is odd. Furthermore, we have $a+1 > c-1$ due to $a+2 < c$. By earlier work, the player on move loses. As Player 2 is on move, Player 1 wins.

Finally, assume $a+2 > c$. Player 1 wins by splitting a two, yielding the game state
\begin{equation}
(a+2, b-1, c).
\end{equation}
Note $a+2$ is odd, $b-1$ is even, and $c$ is even. Moreover, we have $a+2 > c$. By earlier work, the player on move loses. As Player 2 is on move, Player 1 wins.

In either case, we have exhibited a winning strategy for Player 1, and so Player 1 wins by force.
\end{proof}

\subsection{Proof of Theorem \ref{thm:buildup123}}\label{sec:proofthmbuildup123}
\BuildUpOneTwoThree*
\begin{proof}

For $n \equiv 1 \: (4)$ and $n \geq 9$, Player 1 puts down a 3, and then nims down to 6. At this point, the board is either $(O,E,O+1)$ or $(E,E',E+1)$.

\begin{enumerate}[label = (\roman*)]
    \item If Player 2 puts down a 1, the game state is either $(O+1,E,O+1)$ or $(E+1,E',E+1)$.
    \begin{enumerate}[label = (\alph*)]
        \item In the case of $(O+1,E,O+1)$, Player 1 puts down a 2 leaving $(O+1,E+1,O+1)$. 
        \begin{enumerate}[label = (\arabic*)]
           \item If Player 2 puts down a 1 or 2, then Player 1 puts down the other. Player 2 then starts on $(O+2,E+2,O+1)$, which is a win for Player 1. 
           \item If Player 2 puts down a 3, then Player 1 starts on $(O+1,E+1,O+2)$ which is a win for Player 1. 
        \end{enumerate}
        \item In the case of $(E+1,E',E+1)$, Player 1 puts down a 3, leaving $(E+1,E',E+2)$.
        \begin{enumerate}[label = (\arabic*)]
            \item If Player 2 puts down a 1, then player 1 is forced to put down a 1. Player 2 then starts on $(E+3,E',E+2)$, which is a win for Player 1.
            \item If Player 2 puts down a 2, then Player 1 starts on $(E+1,E'+1,E+2)$, which is a win for Player 1.
        \end{enumerate}
    \end{enumerate}
    \item If Player 2 puts down a 2 or 3, then Player 1 puts down the other, forcing Player 2 to put down a 1. Player 1 then starts on $(O+1,E+1,O+2)$ or $(E+1,E'+1,E+2)$, both of which are wins for Player 1. 
\end{enumerate}

For $n \equiv 2 \: (4)$ and $n \geq 6$, Player 2 nims down until there is 6 left. The game state is either $(O,E,O)$ or $(E,E',E)$.
\begin{enumerate}[label = (\roman*)]
    \item If Player 1 puts down a 1, then Player 2 is to move on $(O+1,E,O)$ or $(E+1,E',E)$, with 5 total left.
    \begin{enumerate}[label = (\alph*)]
        \item In the case of $(O+1,E,O)$, Player 2 puts down a 3, leaving $(O+1,E,O+1)$, with a total of 2 left to play.
        \begin{enumerate}[label = (\arabic*)]
            \item If Player 1 puts down a 1, Player 2 is forced to put down a 1. Player 1 then starts on $(O+3,E,O+1)$, which is a win for Player 2.
            \item If Player 1 puts down a 2, Player 2 starts on $(O+1,E+1,O+1)$, which is a win for Player 2.
        \end{enumerate}
        \item In the case of $(E+1,E',E)$, Player 2 puts down a 2,  leaving $(E+1,E'+1,E)$. 
        \begin{enumerate}[label = (\arabic*)]
            \item If Player 1 puts down a 1 or 2, then Player 2 puts the other down. Player 1 then starts on $(E+2,E'+2,E)$, which is a win for Player 2. 
            \item If Player 1 puts down a 3, then Player 2 starts on $(E+1,E'+1,E+1)$, which is a win for Player 2.
        \end{enumerate}
    \end{enumerate}
    \item If Player 1 puts down a 2 or a 3, then Player 2 puts down the other, forcing Player 1 to play a 1. Player 2 then starts on $(O+1,E+1,O+1)$ or $(E+1,E'+1,E+1)$, both of which are wins for Player 2.
\end{enumerate}

For $n \equiv 3 \: (4)$, Player 1 puts down a 2, and then nims down to 1, upon which Player 2 puts down a 1. Player 1 then starts on an ordered triple $(a,b,c)$ with $b$ odd since nimming down preserves the parity of $twos$ and there is a single 2 put down at the start. In all cases where $b$ is odd, the player that starts wins, so Player 1 wins.

For $n = 1$, Player 1 puts down a 1. Player 2 then loses as they cannot move on $(1,0,0)$.

For $n = 2$, if Player 1 puts down a 1, Player 2 must do the same, and so Player 2 wins as Player 1 cannot move on $(2,0,0)$. If Player 2 puts down a 2, Player 2 starts on $(1,1,0)$, and therefore Player 2 wins.

For $n = 4$, Player 1 puts down a 3, forcing Player 2 to put down a 1. Player 1 starts on $(1,0,1)$, a winning position for Player 1.

For $n = 5$, Player 1 puts down a 2.
\begin{enumerate}[label = (\roman*)]
    \item If Player 2 puts down a 3, then Player 1 starts on $(0,1,1)$, a win for Player 1.
    \item If Player 2 puts down a 1 or a 2, then Player 1 puts down the other. Player 2 then starts on $(1,2,0)$, which is a win for Player 1.
\end{enumerate}

\end{proof}

\section{Code}\label{sec:code}

The program we used to brute force solve who has a winning strategy is listed below. It was coded in Jupyter Notebook, but it can be run by any program that can run python. The computation complexity is about $O(\exp(\sqrt{n}))$, with the program taking about 2 hours for $n = 129$ and about 24 hours for $n = 144$.

\begin{lstlisting}[ basicstyle = \ttfamily\small]
import networkx as nx
from matplotlib import pyplot as plt
from networkx.drawing.nx_agraph import graphviz_layout

def is_game_a_loss(current_state):
    """
    Input: A list representing the current state of the game
    Output: Boolean 
    
    Returns True if the game is over (i.e., is a 
        loss for the player next to move)
    Returns False if the game is not over.
    """
    return sum(current_state[1:]) == 0


def combine(current_state):
    """
    Input: A list representing the current state of the game
    Output: A list of lists, where each list is a possible 
        next state of the game
    
    Finds all the potential combine moves in the current 
    state and creates the next states after the combine move
    """
    future_states = []
    for i, val in enumerate(current_state):
        tmp = list(current_state).copy()
        if i >= 3 and val and current_state[i-3]:
            tmp[i] -= 1
            tmp[i-3] -= 1
            tmp[i-1] += 2
            future_states.append(tuple(tmp))
    if len(current_state) > 2:
        if current_state[0] and current_state[2]:
            tmp = list(current_state).copy()
            tmp[0] -= 1
            tmp[2] -= 1
            tmp[1] += 2
            future_states.append(tuple(tmp))
    return set(future_states)  

def split(current_state):
    """
    Input: A list representing the current state of the game
    Output: A list of lists, where each list is a possible next 
        state of the game
    
    Finds all the potential split moves in the current state and creates
    the next states after the split move
    """
    
    if is_game_a_loss(current_state):
        return "An Error Has Occurred"
    future_states = []
    for i, val in enumerate(current_state):
        tmp = list(current_state).copy()
        if i == 0: continue
        if i == 1 and val:
            tmp[1] -= 1
            tmp[0] += 2
            future_states.append(tuple(tmp))
        elif val:
            tmp[i] -= 1
            tmp[i-1] += 1
            tmp[i-2] += 1
            future_states.append(tuple(tmp))
    return set(future_states) 

def nearestSmallerEqFib(n):
    """
    Input: integer n
    Output: tuple, where the first entry is the greatest Fibonacci 
        number smaller than n and the second is the index of 
        that Fibonacci number
    """
     
    # Corner cases
    if (n == 0):
        return "404:BROKEN"
    elif n == 1:
        return (1,n)
        
    # Finds the greatest Fibonacci Number smaller
    # than n.
    f1, f2, f3 = 0, 1, 1
    index = 0
    while (f3 <= n):
        index += 1
        f1 = f2;
        f2 = f3;
        f3 = f1 + f2;
    return (f2, index);

def int_to_zeck(n):
    """
    Input: integer n
    Output: tuple that represents the Zeckendorf decomposition of n
    """
    arr = []
     
    while (n>0):
        f_i, i = nearestSmallerEqFib(n);  
        arr.append(i)
        n = n-f_i
    
    arr2 = [0] * max(arr)
    for i in arr:
        arr2[i-1] += 1    
    return tuple(arr2)

def Fibonacci(n):
    """
    Input: integer n
    Output: integer F_n, the n-th Fibonacci number (where F_2 = 2)
    """
    if n == 0:
        return 1
    # Check if n is 1,2
    # it will return 1
    elif n == 1:
        return 1

    else:
        return (Fibonacci(n-1) + Fibonacci(n-2))
    
def zeck_to_int(zeck_tuple):
    """
    Input: a tuple representing the Zeckendorf decomposition of a number n
    Output: the integer n
    """
    tot = 0
    for i, val in enumerate(zeck_tuple):
        tot += val * Fibonacci(i+1)
        
    return tot
    

def game_solver(number, draw_graph = False, graph_labels = False):  
    """
    Input: an integer n (or a tuple representing a zeckendorf decomposition)
    Output: a tuple of the number, its zeckendorf decomposition, the winner, 
        the number of edges, and vertices
    
    Also can draw the graph of the game
    """
    
    # Checks whether the input is a tuple or an integer, and 
    # converts it to the zeckendorf decomposition if the latter
    if type(number) is tuple:
        initial_state = number
        number = zeck_to_int(number)
    else:
        initial_state = int_to_zeck(number)
        
        
    edges = []      
    current_states = [initial_state]
    calculated_states = []
    
    # Continually loops through the list of possible 
    # states to generate the next states for the game
    while True:
        future_states = []
        for current_state in current_states:
            if current_state in calculated_states:
                continue
            next_states = set()
            if is_game_a_loss(current_state): continue
            next_states = next_states.union(split(current_state))
            next_states = next_states.union(combine(current_state))
            edges += ((current_state, state) for state in next_states)
            future_states += next_states.copy()
            calculated_states += [current_state]
        future_states = set(future_states)
        
        # If there are no states left to be calculated, 
        # we must be at the end node, i.e., game is over
        if len(future_states) == 0:
            break
        current_states = future_states.copy()

    # Initializes the graph of the game
    G = nx.DiGraph()
    for edge in set(edges):
        G.add_edge(edge[0], edge[1])

    # Sets the initial win/loss state for the vertices, all are set to False 
    # (i.e., not calculated) except the end node, which is a loss
    data = dict()
    for node in G.nodes:
        if is_game_a_loss(node):
            data[node] = "L"
        else:
            data[node] = False

    """
    Goes through the game tree in reverse 
    and calculates the state of each vertex
    
    If any vertex leads to one that is a loss, 
        that vertex is automatically a win
    
    If all the child nodes of a vertex are a win, 
        then that vertex is a loss
    """
    not_done = True
    while not_done:
        for node in G.nodes:
            if data[node]: continue
            children_states_L = [data[node] == "L" 
                        for node in G.neighbors(node)]
            children_states_W = [data[node] == "W" 
                        for node in G.neighbors(node)]
            if any(children_states_L):
                data[node] = "W"
            elif all(children_states_W):
                data[node] = "L"

        not_done = not all(data[node] for node in G.nodes)

    # Colors the vertices based on whether they are wins or losses
    color_map = []
    for node in G.nodes:
        if sum(node[1:]) == 0:
            color_map.append("black")
        elif data[node] == "W":
            color_map.append("green")
        elif data[node] == "L":
            color_map.append("red")
        if node == initial_state:
            if data[node] == 'W': winner = 1
            else: winner = 2

    # Draws the graph
    if draw_graph:
        plt.figure(3,figsize=(20,20)) 
        pos=graphviz_layout(G, prog='dot')
        nx.draw_networkx(G, pos, node_color=color_map, with_labels=graph_labels)

    return (number, "".join([str(i) for i in initial_state]), 
        winner, len(G.edges), len(G.nodes))

# Example: should output a 5-tuple (7, '0101', 1, 16, 10) 
#   and a graph representing the game
game_solver(7, True, True)
\end{lstlisting}

%%%%%%%%%%%%%%%%%%%%%%%%%%%%%%%%%%%%%%%%%%%%%%%%%%%%%%%%%%%%%%%%%%%%%%%%%%%%%%%%%%%%%%%%%%%%%%%%%%%%%%%%%%%%%%%%%%%%%%%%
%%%%%%%%%%%%%%%%%%%%%%%%%%%%%%%%%%%%%%%%%%%%%%%%%%%%%%%%%%%%%%%%%%%%%%%%%%%%%%%%%%%%%%%%%%%%%%%%%%%%%%%%%%%%%%%%%%%%%%%%
%%%%%%%%%%%%%%%%%%%%%%%%%%%%%%%%%%%%%%%%%%%%%%%%%%%%%%%%%%%%%%%%%%%%%%%%%%%%%%%%%%%%%%%%%%%%%%%%%%%%%%%%%%%%%%%%%%%%%%%%

\end{document}